\documentclass{amsart}
\usepackage[margin=1.5in]{geometry}

\usepackage{amsmath,amsthm, latexsym, amssymb,mathrsfs,  tikz-cd}
\usepackage{stackengine,scalerel, quiver, MnSymbol}
\usepackage[colorlinks]{hyperref}
\definecolor{dark-red}{rgb}{0.6,0,0}
\definecolor{dark-green}{rgb}{0,0.4,0}
\definecolor{medium-blue}{rgb}{0,0,0.5}
\hypersetup{
    colorlinks, linkcolor={dark-red},
    citecolor={dark-green}, urlcolor={medium-blue}
}

\usepackage{relsize}
\usepackage[bbgreekl]{mathbbol}
\usepackage{amsfonts}
\DeclareSymbolFontAlphabet{\mathbb}{AMSb} 
\DeclareSymbolFontAlphabet{\mathbbl}{bbold}

\newcommand{\mr}[1]{\mathrm{#1}}

\newcommand{\BC}{\mathrm{BC}}

\newcommand{\Spa}{\mathrm{Spa}}

\newcommand{\crys}{\mr{crys}}

\numberwithin{equation}{subsection}
\numberwithin{equation}{subsubsection}

\theoremstyle{plain}

\newtheorem*{theorem*}{Theorem}
\newtheorem*{conjecture*}{Conjecture}

\theoremstyle{definition}

\newtheorem{example}[subsubsection]{Example}

\title{The relativistic $p$-adic sunscreen conjecture}

\author{Sean Howe}

\begin{document}

\begin{abstract} We formulate a conjecture about intersections between the Banach-Colmez space $\BC(1/2)$ and germs of smooth rigid analytic curves at the origin in $\mathbb{A}^2_{\mathbb{C}_p}$.  
\end{abstract}

\maketitle

\section{The conjecture}

Let $p$ be a prime number, let $\mathbb{Q}_p$ be the $p$-adic numbers, and let $\mathbb{C}_p$ be the completion of an algebraic closure of $\mathbb{Q}_p$ for the unique extension of the $p$-adic absolute value. 

\subsection{$\BC(1/2)$ as $p$-adic sunscreen}
The Banach--Colmez space $\BC(1/2)$ admits (over $\mathbb{C}_p$) a natural embedding into the rigid analytic plane $\mathbb{A}^2_{\mathbb{C}_p}$. It can be described in several ways: if $\BC(1/2)$ is viewed as the global sections of the rank two vector bundle $\mathcal{O}(1/2)$ on the Fargues--Fontaine curve $X_{\mathbb{C}_p^\flat}$, then the embedding is by evaluation at the canonical point $\infty:\; \Spa \mathbb{C}_p \hookrightarrow X_{\mathbb{C}_p^\flat}$. If $\BC(1/2)$ is viewed as the inverse limit under multiplication by $p$ for the points of the generic fiber of the Lubin-Tate formal group for $\mathbb{Q}_{p^2}$, then the embedding is by the quasi-logarithm. Most explicitly, if we consider Fontaine's crystalline period ring $B^+_\crys$ with the natural surjection $\theta: B^+_\crys \twoheadrightarrow \mathbb{C}_p$ and Frobenius endomorphism $\varphi: B^+_\crys \rightarrow B^+_\crys$, then $\BC(1/2)=B_\crys^{+, \varphi^2=p}$ and the embedding into $\mathbb{C}_p^2$ sends $b \in B_\crys^{+, \varphi^2=p}$ to $\left(\theta(b), \theta(\varphi(b))\right)$. 

It will suffice below to consider just the $\mathbb{C}_p$-points, so we view $\BC(1/2)$ as a subset of $\mathbb{C}_p^2$. It is an infinite dimensional $\mathbb{Q}_p$-Banach subspace with the following remarkable property:\\

\noindent{\bf The $p$-adic sunscreen property.} For any $\mathbb{C}_p$-line $\ell \subseteq \mathbb{C}_p^2$, $\ell \cap \BC(1/2)$ is a $\mathbb{Q}_p$-plane, i.e. a translate of a 2-dimensional $\mathbb{Q}_p$-subspace of the $\mathbb{Q}_p$-Banach space $\BC(1/2)$. \\

In particular, if you stepped outside on a two-dimensional $p$-adic planet coated in a layer of $\BC(1/2)$ then, although you would barely notice it on you as more than a thin 1-dimensional $\mathbb{C}_p$-sheet, you would be well-protected from arriving UV-rays: from any angle of entry, the rays are blocked not just by one particle but by an entire profinite set's worth!  

Viewing $\BC(1/2)$ as the universal cover of the Lubin-Tate formal group, the $p$-adic sunscreen property is a reformulation of the surjectivity of the Gross-Hopkins period map to $\mathbb{P}^1$ combined with the description of the $p$-adic comparison theorem via the quasi-logarithm as in \cite{ScholzeWeinstein.ModuliOfpDivisibleGroups}. In terms of the Fargues--Fontaine curve, it is the statement that any modification of $\mathcal{O}(1/2)$ by a $1$-dimensional subspace of $\infty^*\mathcal{O}(1/2)$ is the trivial rank two vector bundle, combined with the computation of $H^\bullet(X_{\mathbb{C}_p^\flat}, \mathcal{O})$ (i.e. the fundamental exact sequence). 

\subsection{Relativistic accommodations}
This is a mathematical solution to $p$-adic sunburn, but, unfortunately, not an engineering solution. The reason being, of course, that the theory of general relativity predicts that the gravity well of the $p$-adic planet you are standing on will slightly curve the incoming rays of light. Thus, it is not enough to block every line with a profinite set; $\BC(1/2)$ will only work as sunscreen if it can also obstruct curves. 

Thus, before we bring our $p$-adic sunscreen to market, we must resolve:
\begin{conjecture*}[The relativistic $p$-adic sunscreen conjecture for $\BC(1/2)$]
Let $F(x,y) \in \mathbb{C}_p[[x,y]]$ be a power series with constant term $0$ and non-vanishing gradient (i.e. $F(x,y)=0+ ax + by + \ldots$ with $(a,b)\neq (0,0)$). Then, in any disk $D=\{(x,y) \in \mathbb{C}_p^2, |x| \leq \delta, |y| \leq \delta)$ where $F$ converges, the set
$\{ (x,y) \in \BC(1/2) \cap D\, |\, F(x,y)=0 \}$
is a profinite set of cardinality $2^\mathbb{N}$. 
\end{conjecture*}

The $p$-adic sunscreen property implies this conjecture for linear $F$. Explicitly, points in the intersection can be interpreted as solutions $b\in B^{+,\varphi^2=p}_\crys$ to $F(\theta(b), \theta(\varphi(b)))=0$. 

\begin{example}
    To our knowledge, this conjecture is open for any non-linear convergent $F$. In particular, one can take, $F(x,y)=y^2-x$, whose vanishing set is the parabola $y^2=x$ and, for this $F$, it is not even clear whether there is \emph{any} other solution except $(0,0)$. 
\end{example}

\noindent{\bf Bounty:} The $p$-adic sunscreen story was inspired by Falconer's work on projections in fractal geometry. I will send you a digital sundial if you resolve this conjecture for $y^2=x$.

\subsection{Heuristic}
As is typical with general relativity, in order to understand why the conjecture might be true one needs to introduce new ideas in differential geometry. Recent computations (cf., e.g., \cite{Howe.InscriptionTwistorsAndPAdicPeriods}) show that it may be reasonable to view $\BC(1/2)$ as some type of Banach--Colmez submanifold of the rigid analytic variety $\mathbb{C}_p^2$; because it is a $\mathbb{Q}_p$-vector space object, its tangent space at any point should be naturally identified with itself, i.e. with $\BC(1/2) \subseteq \mathbb{C}_p^2$. On the other hand, the vanishing set of $F(x,y)$ in any small disk is a smooth rigid analytic curve $C$ containing $(0,0)$, so has a natural tangent line $T_{C,(0,0)} \subseteq \mathbb{C}_p^2$. 

The $p$-adic sunscreen property implies that $T_{C,(0,0)} + \BC(1/2) = \mathbb{C}_p^2$, and thus the intersection of $C$ and $\BC(1/2)$ at $(0,0)$ should be \emph{transverse}. Transverse intersections of Banach--Colmez manifolds should be manifolds, and, moreover the putative tangent space of the intersection, $T_{C,(0,0)} \cap \BC(1/2),$ is a 2-dimensional $\mathbb{Q}_p$-vector space. Thus we expect the intersection should look, locally around 0, like a $2$-dimensional $\mathbb{Q}_p$-analytic manifold; instead of trying to make the $p$-adic manifold structure precise, we are content for now to conjecture just the topological property that one finds a profinite set of the same cardinality. 

\subsection{Immediate generalizations}
For any $n \geq 2$, the same considerations can be made for $\BC(1/n)\subseteq \mathbb{C}_p^n$ by intersection with smooth hypersurface germs at the origin. More generally, one could formulate a version for any Banach--Colmez space arising as the global sections of a vector bundle on the Fargues--Fontaine curve with slopes between $0$ and $1$ (or, equivalently, as the universal cover of a $p$-divisible group) intersecting smooth rigid analytic varieties of appropriate dimensions. However the local structure predicted is more complicated outside of the $\BC(1/n)$ case and, in some cases, one has to forbid certain tangent directions to ensure transversality. For example, one can intersect $\BC(1/2)^2 \subseteq (\mathbb{C}_p^2)^2$ with a surface $S$, but the intersection will not be transverse if the tangent space $T_{S,\vec{0}}$ is equal to the first copy of $\mathbb{C}_p^2$: in this case, the quotient $\mathbb{C}_p^4/(\BC(1/2)^2 + T_{S,\vec{0}})\cong\mathbb{C}_p^2/\BC(1/2)$ can be identified with $\BC(-1/2)$, which is the quotient of $\mathbb{C}_p$ by a two-dimensional $\mathbb{Q}_p$-subspace.

\subsection{Why?}
One can formulate analogs of these questions for intersections between rigid analytic varieties and other perfectoid spaces or diamonds in the place of $\BC(1/2)$. In particular, in \cite{Howe.InscriptionTwistorsAndPAdicPeriods} we have introduced a theory of ``inscribed $v$-sheaves" where these differential computations are valid mathematics and, for example, in \cite[\S6]{HoweKlevdal.AdmissiblePairsAndpAdicHodgeStructuresIIIVariationAndUnlikelyIntersection}  we use this notion of transversality to formulate a $p$-adic Ax-Schanuel conjecture for infinite level local Shimura varieties. However, this theory is based on nilpotents, thus it is completely unclear how the differential computations reflect the local geometry of actual points. The relativistic $p$-adic sunscreen conjecture attempts to clarify this relation in the simplest non-trivial case.

\subsection{Acknowledgements} We thank Pierre Colmez, Peter Wear, and Jared Weinstein for helpful discussions. We have stated this conjecture in some communciations previously, but we gave it this fun name and backstory at the 2026 \emph{Geometrization of the Local Langlands Correspondence} conference at CIRM; we thank the organizers and staff (especially Joaquin and that one guy who none of us will ever forget) for the opportunity and lovely conference! 
\bibliographystyle{plain}
\bibliography{references, preprints}

\end{document}